\documentclass{article}
\usepackage{amsmath,amssymb}
\begin{document}
\title{A succinct method for investigating the sums of infinite series through differential formulae\footnote{Delivered to the St. Petersburg Academy
on March 13, 1780. Originally published as {\em Methodus succincta summas serierum infinitarum per formulas differentiales investigandi}, M\'emoires de l'Acad\'emie Imp\'eriale des Sciences de St.-P\'etersbourg \textbf{5} (1815), 45--56, and reprinted in \emph{Leonhard Euler, Opera Omnia}, Series 1:
Opera mathematica,
Volume 16.
A copy of the original text is available
electronically at the Euler Archive, at http://www.eulerarchive.org. This paper
is E746 in the Enestr\"om index.}}
\author{Leonhard Euler\footnote{Date of translation: May 5, 2007. Translated from the Latin by Jordan Bell, 4th year undergraduate in Honours Mathematics, School of Mathematics and Statistics, Carleton University, Ottawa, Ontario, Canada. Email: jbell3@connect.carleton.ca}}
\date{}

\maketitle

\S 1. Even though I have already dealt with this matter several times, 
most of what has been found for conveniently expressing sums has been dispersed in many different papers, and
moreover been done somewhat ambiguously;
because of this, I will relate this succinct method,
through which the sum of any series can be easily calculated, without ambiguity,
by a very simple form which will be worked out.

\S 2. Thus let $X$ be some function of $x$, and let $X',X'',X'''$, etc. arise from it by respectively writing $x+1,x+2,x+3$, etc. in place of $x$. 
Then here the letters $X,X',X'',X'''$, etc. designate to me the terms of some series corresponding to the indices $x,x+1,x+2,x+3$, etc.
I will consider two cases of infinite series in these forms; in the first, the
terms all appear with the same sign $+$, so that the series which is to be summed would be:
\[
X+X'+X''+X'''+\textrm{etc.}
\]
Then in the other case, these same terms proceed with alternating signs, so that the series which is to be summed would be $X-X'+X''-X'''+$ etc. Accordingly, I will now deal with these two cases, one then the other.

\begin{center}
Case 1.
Summation of the infinite series
\[
S=X+X'+X''+X'''+\textrm{etc.}
\]
\end{center}

\S 3. Let $S'$ denote the sum of this series with the first term truncated,
that is, so that $S'=X'+X''+X'''+$ etc. Since $S$ is a certain function of $x$, which indeed we are investigating
here principally, $S'$ will likewise be a function of $x+1$. It is thus evident that $S-S'=X$. Now, since
\[
S'=S+\partial S+\frac{1}{2}\partial \partial S+\textrm{etc.},
\]
where I have chosen for brevity to suppress the denominators containing powers of the element $\partial x$, it is seen at once that our series will assume this form:
\[
0=X+\partial S+\frac{1}{2}\partial \partial S+\frac{1}{6}\partial^3 S+\frac{1}{24}\partial^4 S+\textrm{etc.}
\]

\S 4. Therefore if that series converges strongly, then it will approximately be
$\partial S=-X$, and so $S=-\int X\partial x$; but this integral should be determined by a constant that vanishes when $x$ is taken infinitely large, since
the infinitesimal terms can be considered as nothing, as the latter series has no finite value. Having considered the sum approximately, for the true sum we set
\[
S=-\int X\partial x-\alpha X-\beta \partial X-\gamma \partial \partial X-\textrm{etc.}
\]
and then it will be
\[
\partial S=-X-\alpha \partial X-\beta \partial \partial X-\gamma \partial^3 X-\textrm{etc.}
\]
But if the values arising here are substituted for each of the differentials of $S$, the following equation will be obtained:
\[
\left.
\begin{array}{llllll}
+X&-\alpha \partial X&-\beta \partial \partial X&-\gamma \partial^3 X&-\delta \partial^4 X&-\textrm{etc.}\\
-X&-\frac{1}{2}&-\frac{1}{2}\alpha&-\frac{1}{2}\beta&-\frac{1}{2}\gamma
&-\textrm{etc.}\\
&&-\frac{1}{6}&-\frac{1}{6}\alpha&-\frac{1}{6}\beta&-\textrm{etc.}\\
&&&-\frac{1}{24}&-\frac{1}{24}\alpha&-\textrm{etc.}\\
&&&&-\frac{1}{120}&-\textrm{etc.}\\
\end{array} \right\}=0
\]
and now the unknown coefficients $\alpha,\beta,\gamma$, etc. can be defined by the following equalities:
\[
\alpha+\frac{1}{2}=0; \quad \beta+\frac{1}{2}\alpha+\frac{1}{6}=0; \quad \gamma+\frac{1}{2}\beta+\frac{1}{6}\alpha+\frac{1}{24}=0; \quad \textrm{etc.},
\]
from which we shall have $\alpha=-\frac{1}{2}$; $\beta=\frac{1}{12}$; $\gamma=0$; etc.

\S 5. However, this way of evaluating the letters $\alpha,\beta,\gamma$, etc. takes a fair bit of work, and
also, there is not an obvious rule for proceeding further; thus I will now inquire into the values of these letters in a straightforward way.
Namely I will consider a series formed with the same coefficients as the preceding.
Let
$V=1+\alpha z+\beta z^2+\gamma z^3+\delta z^4+$ etc. Indeed it is evident that should the sum of this series $V$ be brought
to a finite form and each of the powers of $z$ expanded, the same series will necessarily hold. 
Then having worked this out, the values of the letters $\alpha,\beta,\gamma,\delta$, etc. will 
freely reveal themselves.

\S 6. Therefore from the relations which exist between the letters $\alpha,\beta,\gamma,\delta$, etc., which were
written out above in \S 4, we get the following:
\[
\begin{array}{rclllllll}
V&=&1&+\alpha z&+\beta z^2&+\gamma z^3&+\delta z^4&+\epsilon z^5&+\textrm{etc.}\\
\frac{1}{2}z V&=&&+\frac{1}{2}&+\frac{1}{2}\alpha&+\frac{1}{2}\beta&+\frac{1}{2}\gamma&+\frac{1}{2}\delta&+\textrm{etc.}\\
\frac{1}{6}zz V&=&&&+\frac{1}{6}&+\frac{1}{6}\alpha&+\frac{1}{6}\beta&+\frac{1}{6}\gamma&+\textrm{etc.}\\
\frac{1}{24}z^3 V&=&&&&+\frac{1}{24}&+\frac{1}{24}\alpha&+\frac{1}{24}\beta&+\textrm{etc.}\\
\frac{1}{120}z^4 V&=&&&&&+\frac{1}{120}&+\frac{1}{120}\alpha&+\textrm{etc.}\\
\frac{1}{720}z^5 V&=&&&&&&+\frac{1}{720}&+\textrm{etc.}\\
&&\textrm{etc.}&&&&&&
\end{array}
\]
Of course here all the terms, aside from the first, cancel each other out; it will therefore be
\[
V(1+\frac{1}{2}z+\frac{1}{6}z^3+\frac{1}{24}z^3+\frac{1}{120}z^4+\frac{1}{720}z^5+\textrm{etc.})=1.
\]

\S 7. Thus since $e^z=1+z+\frac{1}{2}z^2+\frac{1}{6}z^3+\frac{1}{24}z^4+$ etc. it will be $\frac{V(e^z-1)}{z}=1$, and thus $V=\frac{z}{e^z-1}$; to make it easier to convert this expression again into a series, let us put $z=2t$, so that it would be $V=\frac{2t}{e^{2t}-1}$, and therefore $V+t=t\cdot \frac{e^{2t}+1}{e^{2t}-1}$. Now setting $\frac{e^{2t}+1}{e^{2t}-1}=u$ it follows that $V=tu-t$. Therefore since $u=\frac{e^t-e^{-t}}{e^t-e^{-t}}$, by expanding the exponentials it will be
\[
u=\frac{1+\frac{1}{2}t^2+\frac{1}{24}t^4+\frac{1}{720}t^6+\textrm{etc.}}{t+\frac{1}{6}t^3+\frac{1}{120}t^5+\frac{1}{5040}t^7+\textrm{etc.}},
\]
where in the numerator only even powers, and indeed in the denominator only odd powers, appear. Moreover, it is apparent that by making $t$ small enough it will
make $u=\frac{1}{t}$, and indeed it also follows for the terms to progress as the powers $t,t^3,t^5$, etc.

\S 8. Therefore since we have put $u=\frac{e^{2t}+1}{e^{2t}-1}$, it will be
$e^{2t}=\frac{u+1}{u-1}$, and then $2t=l\frac{u+1}{u-1}$. Therefore by differentiating here it will be $\partial t=-\frac{\partial u}{uu-1}$, from which it is concluded $\frac{\partial u}{\partial t}+uu-1=0$. Moreover, because we know
that the first term of the series by which $u$ is expressed is $\frac{1}{t}$ and that the exponents of the successive powers increase by two, it may be set:
\[
u=\frac{1}{t}+2At-2Bt^3+2Ct^5-2Dt^7+\textrm{etc.}
\] 
and the substitution happens in the following way:
\[
\begin{array}{rclllllll}
\frac{\partial u}{\partial t}&=&-\frac{1}{tt}&+2A&-6Btt&+10Ct^4&-14Dt^6&+18Et^8&-\textrm{etc.}\\
uu&=&+\frac{1}{tt}&+4A&-4B&+4C&-4D&+4E&-\textrm{etc.}\\
&&&&+4AA&-8AB&+8AC&-8AD&+\textrm{etc.}\\
&&&&&&+4BB&-8BC&+\textrm{etc.}\\
-1&=&-1&&&&&&
\end{array}
\]
where the first terms cancel each other out, and then indeed the remaining obtain the following determinations:
\begin{eqnarray*}
6A=1&\textrm{therefore}&A=\frac{2}{3}\cdot \frac{1}{4}=\frac{1}{6},\\
10B=4AA&&B=\frac{2}{5}AA=\frac{1}{90},\\
14C=8AB&&C=\frac{2}{7}\cdot 2AB=\frac{1}{945},\\
18D=8AC+4BB&&D=\frac{2}{9}(2AC+BB)=\frac{1}{9450},\\
22E=8(AD+BC),&&E=\frac{2}{11}(2AD+2BC)=\frac{1}{93555},\\
&\textrm{etc.}&
\end{eqnarray*}

\S 9. Therefore these letters $A,B,C,D$, etc. are entirely the same as those which I previously used for expressing the sums of the reciprocals of powers, 
which were accordingly found to be:
\begin{eqnarray*}
1+\frac{1}{4}+\frac{1}{9}+\frac{1}{16}+\frac{1}{25}+\textrm{etc.}&=&A\pi^2,\\
1+\frac{1}{4^2}+\frac{1}{9^2}+\frac{1}{16^2}+\frac{1}{25^2}+\textrm{etc.}&=&B\pi^4,\\
1+\frac{1}{4^3}+\frac{1}{9^3}+\frac{1}{16^3}+\frac{1}{25^3}+\textrm{etc.}&=&C\pi^6,\\
\textrm{etc.}&&
\end{eqnarray*}
and I have already worked out these values all the way to the thirty fourth by very painstaking calculation.

\S 10. Since we have supposed that
\[
u=\frac{1}{t}+2At-2Bt^3+2Ct^5-\textrm{etc.}
\]
and because $V=tu-t$, it will be
\[
V=1-t+2At^2-2Bt^4+2Ct^6-2Dt^8+\textrm{etc.}.
\]
Here there is not much else for us to do, except that as $t$ is written in place of $\frac{1}{2}z$, it follows that
\[
V=1-\frac{z}{2}+\frac{Azz}{2}-\frac{Bz^4}{8}+\frac{Cz^6}{32}-\frac{Dz^8}{128}+\textrm{etc.}
\]
On the other hand since we have
\[
V=1+\alpha z+\beta z^2+\gamma z^3+\textrm{etc.},
\]
by collecting terms together we will find that $\alpha=\frac{1}{2}$; $\beta=\frac{1}{2}A$; $\gamma=0$; $\delta=-\frac{1}{8}B$; $\epsilon=0$; $\zeta=\frac{1}{32}C$; $\eta=0$; etc.

\S 11. With the values of these letters now found, the sum of the given series
\[
S=X+X'+X''+X'''+\textrm{etc.}
\]
can be expressed in the following way:
\[
S=-\int X \partial x + \frac{1}{2}X-\frac{1}{2}A \partial X+\frac{1}{8}B\partial^3 X-\frac{1}{32}C\partial^5 X+\frac{1}{128}D\partial^7 X-\frac{1}{512}E\partial^9 X+\textrm{etc.},
\]
where the integral $\int X \partial x$ should be taken such that it vanishes by putting $x=\infty$; it is then clear
that if the constant that is to be added has to be infinite, then likewise the sum of this series will be infinite.

\S 12. Let us consider the example in which $X=\frac{1}{x^n}$, so that the sum of this series will be sought:
\[
S=\frac{1}{x^n}+\frac{1}{(x+1)^n}+\frac{1}{(x+2)^n}+\frac{1}{(x+3)^n}+\textrm{etc.}
\]
Here it will therefore be $\int X \partial x=-\frac{1}{(n-1)x^{n-1}}$. For this form to vanish by putting $x=\infty$ it is necessary that the exponent $n$ be greater than unity; for otherwise, were it $n=1$ or $n<1$, the sum of
the series would certainly be infinitely large. 
Now, it will indeed be $\partial X=-\frac{n}{x^{n+1}}$, and then $\partial^3 X=-\frac{n(n+1)(n+2)}{x^{n+3}}$; 
$\partial^5 X=-\frac{n\cdots (n+4)}{x^{n+5}}$; etc., and by substituting in these values, the sum that is being sought out will be:
\[
S=\frac{1}{(n-1)x^{n-1}}+\frac{1}{2x^n}+\frac{A}{2}\cdot \frac{n}{x^{n+1}}-\frac{B}{8}\cdot \frac{n(n+1)(n+2)}{x^{n+3}}+
\frac{C}{32}\cdot \frac{n\cdots (n+4)}{x^{n+5}}-\textrm{etc.}
\]
This series converges more strongly the larger the number that is taken for $x$, and this is in addition to the fact that
the letters $A,B,C$, etc. themselves constitute a rapidly convergent progression.

\S 13. Therefore if the terms starting at unity are gathered together $1+\frac{1}{2^n}+\frac{1}{3^n}+\frac{1}{4^n}+\ldots+\frac{1}{(x-1)^n}$, and their sum is called $\Delta$, the sum of the same series continued to infinity will be
$\Delta+S$. It was in this way that I had formerly computed the sums of such infinite series for
the exponent $n$ with each of the values $2,3,4,5$, etc. to many decimal places, by taking namely $x=10$. With this
done, the calculation could be carried out fairly easily.

\begin{center}
Case 2.
Summation of the infinite series
\[
S=X-X'+X''-X'''+X^{IV}-\textrm{etc.}
\]
\end{center}

\S 14. Thus if the index $x$ is increased by unity, we will have
$S'=X'-X''+X'''-X^{IV}+$ etc. This equation is added to the preceding one, and the finite equation $S+S'=X$ is obtained.
From this we will have, by differential formulae,
\[
X=2S+\partial S+\frac{1}{2}\partial \partial S+\frac{1}{6}\partial^3 S+\frac{1}{24}\partial^4 S+\textrm{etc.},
\]
where if the differentials are neglected, it will be $S=\frac{1}{2}X$, which will therefore
be the first term of the series we are searching for. We therefore set
\[
S=\frac{1}{2}X+\alpha \partial X+\beta \partial \partial X+\gamma \partial^3 X+\textrm{etc.}
\]
and having done the substitution it becomes:
\[
\begin{array}{rcllllll}
2S&=&X&+2\alpha \partial X&+2\beta \partial \partial X&+2\gamma \partial^3 X&+2\delta \partial^4 X&+\textrm{etc.}\\
\partial S&=&&+\frac{1}{2}&+\alpha&+\beta&+\gamma&+\textrm{etc.}\\
\frac{1}{2}\partial \partial S&=&&&+\frac{1}{4}&+\frac{1}{2}\alpha&+\frac{1}{2}\beta&+\textrm{etc.}\\
\frac{1}{6}\partial^3 S&=&&&&+\frac{1}{12}&+\frac{1}{6}\alpha&+\textrm{etc.}\\
\frac{1}{24}\partial^4 S&=&&&&&+\frac{1}{48}&+\textrm{etc.}\\
\textrm{etc.}&&&&&&\textrm{etc.}&
\end{array}
\]
where the entire expression is equal to the single term $X$. 

\S 15. Thus with each of the vertical columns reduced to nothing, the following equalities arise:
\[
2\alpha+\frac{1}{2}=0; \, 2\beta+\alpha+\frac{1}{4}=0; \, 2\gamma+\beta+\frac{1}{2}\alpha+\frac{1}{12}=0;
\, 2\delta+\gamma+\frac{1}{2}\beta+\frac{1}{6}\alpha+\frac{1}{24}=0; \, \textrm{etc.}
\]
Here like before the letters take these determinations:
\[
\alpha=-\frac{1}{4}; \, \beta=0; \, \gamma=\frac{1}{4}; \, \delta=0; \, \textrm{etc.}
\]

\S 16. To help us investigate these values, we shall consider the series:
\[
V=\frac{1}{2}+\alpha z+\beta z^2+\gamma z^3+\textrm{etc.}
\]
where now the sum $V$ needs to be searched for. We then derive the following series:
\[
\begin{array}{rclllllll}
2V&=&1&+2\alpha z&+2\beta zz&+2\gamma z^3&+2\delta z^4&+2\epsilon z^5&+\textrm{etc.}\\
Vz&=&&+\frac{1}{2}z&+\alpha zz&+\beta z^3&+\gamma z^4&+\delta z^5&+\textrm{etc.}\\
\frac{1}{2}Vzz&=&&&+\frac{1}{4}&+\frac{1}{2}\alpha&+\frac{1}{2}\beta&+\frac{1}{2}\gamma&+\textrm{etc.}\\
\frac{1}{6}Vz^3&=&&&&+\frac{1}{12}&+\frac{1}{6}\alpha&+\frac{1}{6}\beta&+\textrm{etc.}\\
\frac{1}{24}Vz^4&=&&&&&+\frac{1}{48}&+\frac{1}{24}\alpha&+\textrm{etc.}\\
&&\textrm{etc.}&&&&&&
\end{array}
\]
Thus the sum of these series, by the equalities related before, will be $=1$, and so we will have this equation:
\[
V(2+z+\frac{1}{2}z^2+\frac{1}{6}z^3+\frac{1}{24}z^4+\textrm{etc.})=1.
\]
Now, since 
\[
e^z=1+z+\frac{1}{2}z^2+\frac{1}{6}z^3+\textrm{etc.}
\]
it will clearly be $V(1+e^z)=1$, or $V=\frac{1}{1+e^z}$, from which
it is $2V-1=\frac{1-e^z}{1+e^z}$.

\S 17. Thus if it is put as before $\frac{e^z-1}{e^z+1}=u$, so that it would
be $2V=1-u$, and then again $z=2t$ so that $u=\frac{e^t-e^{-t}}{e^t+e^{-t}}$,
by expanding this it will be $u=\frac{t+\frac{1}{6}t^3+\frac{1}{120}t^5+\frac{1}{5040}t^7+\textrm{etc.}}{1+\frac{1}{2}t^2+\frac{1}{24}t^4+\frac{1}{720}t^6+\textrm{etc.}}$. It is apparent from this series that the first term of the value
$u$ expresses will be $t$, and indeed that the following will proceed in odd powers of $t$.

\S 18. Since it is $u=\frac{e^{2t}-1}{e^{2t}+1}$, it will be $e^{2t}=\frac{1+u}{1-u}$, and then
$2t=l\frac{1+u}{1-u}$. By differentiating, it will be $\partial t=\frac{\partial u}{1-uu}$, so that,
$\frac{\partial u}{\partial t}+uu-1=0$, which is the same equation found in the first case. However
the same series for $u$ does not occur. Because the first term of this series should be $=t$, the series will take this type of form:
\[
u=t-\mathfrak{A}t^3+\mathfrak{B}t^5-\mathfrak{C}t^7+\mathfrak{D}t^9-\mathfrak{E}t^{11}+\textrm{etc.}
\]
and by substitution it will become:
\[
\begin{array}{rclllllll}
\frac{\partial u}{\partial t}&=&1&-3\mathfrak{A}tt&+5\mathfrak{B}t^4&-7\mathfrak{C}t^6&+9\mathfrak{D}t^8&-11\mathfrak{E}t^{10}&+\textrm{etc.}\\
uu&=&&+1&-2\mathfrak{A}&+2\mathfrak{B}&-2\mathfrak{C}&+2\mathfrak{D}&-\textrm{etc.}\\
&&&&&\mathfrak{A}^2&-2\mathfrak{A}\mathfrak{B}&+2\mathfrak{A}\mathfrak{C}&-\textrm{etc.}\\
&&&&&&&+\mathfrak{B}^2&-\textrm{etc.}\\
-1&=&-1&&&&&&
\end{array}
\]
From this, the following determinations then follow:
\begin{eqnarray*}
3\mathfrak{A}=1&\textrm{and so}&\mathfrak{A}=\frac{1}{3},\\
5\mathfrak{B}=2\mathfrak{A}&\textrm{and so}&\mathfrak{B}=\frac{2}{5}\mathfrak{A}=\frac{2}{15},\\
7\mathfrak{C}=2\mathfrak{B}+\mathfrak{A}^2,&\textrm{hence}&\mathfrak{C}=\frac{2}{7}\mathfrak{B}+\frac{1}{7}\mathfrak{A}^2=\frac{17}{315},\\
9\mathfrak{D}=2\mathfrak{C}+2\mathfrak{A}\mathfrak{B}&\textrm{thus}&\mathfrak{D}=\frac{2}{9}\mathfrak{C}+\frac{2}{9}\mathfrak{A}\mathfrak{B}=\frac{62}{2835},\\
\textrm{etc.}&&\textrm{etc.}
\end{eqnarray*}

\S 19. Therefore since it is $V=\frac{1}{2}-\frac{1}{2}u$, if in place of $t$ we replace $\frac{z}{2}$, then we obtain this series for $V$:
\[
V=\frac{1}{2}-\frac{1}{4}z+\frac{1}{16}\mathfrak{A}z^3-\frac{1}{64}\mathfrak{B}z^5+\frac{1}{256}\mathfrak{C}z^7-\frac{1}{1024}\mathfrak{D}z^9+\textrm{etc.}
\]
Now since we had set
\[
V=\frac{1}{2}+\alpha z+\beta z^2+\gamma z^3+\delta z^4+\textrm{etc.}
\]
then we collect the values of the letters $\alpha,\beta,\gamma,\delta$, etc., which will therefore
be $\alpha=-\frac{1}{4}$; $\beta=0$; $\gamma=\frac{1}{16}\mathfrak{A}$; $\delta=0$; $\epsilon=-\frac{1}{64}\mathfrak{B}$;
$\zeta=0$; $\eta=\frac{1}{256}\mathfrak{C}$; $\theta=0$; etc. Consequently, the sum that is being sought out will be:
\[
S=\frac{1}{2}X-\frac{1}{4}\partial X+\frac{1}{16}\mathfrak{A}\partial^3 X+\frac{1}{64}\mathfrak{B}\partial^5 X+\frac{1}{256}\mathfrak{C}\partial^7 X-\textrm{etc.}
\]

\S 20. Let us now compare these coefficients to those in the preceding case, which we had obtained by similar differentials. They were $\frac{A}{2},\frac{B}{8},\frac{C}{32}$, etc., and we will discover a remarkable relation between the two, as can be seen from this scheme:
\[
\begin{array}{r|l}
\partial X&\frac{1}{4}:\frac{A}{2}=3=2^2-1,\\
\partial^3 X&\frac{\mathfrak{A}}{16}:\frac{B}{8}=15=2^4-1,\\
\partial^5 X&\frac{\mathfrak{B}}{64}:\frac{C}{32}=63=2^6-1,\\
\partial^7 X&\frac{\mathfrak{C}}{256}:\frac{D}{128}=255=2^8-1,\\
\partial^9 X&\frac{\mathfrak{D}}{1024}:\frac{E}{512}=1023=2^{10}-1,\\
\textrm{etc.}&\textrm{etc.}
\end{array}
\]

\S 21. Thus also in this case, the sum that is being sought out can be conveniently expressed by the very same well known numbers $A,B,C,D$, etc. in the following way:
\begin{eqnarray*}
S&=&\frac{1}{2}X-(2^2-1)\frac{A}{2}\cdot \partial X+(2^4-1)\frac{B}{8}\cdot \partial^3 X-(2^6-1)\frac{C}{32}\cdot \partial^5 X\\
&&+(2^8-1)\frac{D}{128}\cdot \partial^7 X-(2^{10}-1)\frac{E}{512}\cdot \partial^9 X+\textrm{etc.},
\end{eqnarray*}
and this series can be continued as far as we like.

\end{document}